\documentclass{amsart}
\usepackage{latexsym,amssymb,amsthm,amsmath}

\usepackage{amssymb}
\usepackage{amsmath}

\theoremstyle{plain}
\newtheorem{theorem}{Theorem}[section]
\newtheorem*{Theorem B}{Theorem B}
\newtheorem*{Theorem A}{Theorem A}

\newtheorem{definition}{Definition}[section]

\numberwithin{equation}{section}

\theoremstyle{remark}

\begin{document}
\title[Some remarks on Yamabe solitons ]
{Some remarks on Yamabe solitons}
\author[D. Chakraborty, Y. C. Mandal and S. K. Hui$^*$]{Debabrata Chakraborty, Yadab Chandra Mandal and Shyamal Kumar Hui$^*$}
\keywords{Yamabe solitons, Riemannian manifold, Killing vector field.\\ $^*$ Corresponding author}
\begin{abstract}
In this paper we have obtained evolution of some geometric quantities on a compact Riemannian manifold $M^n$ when the metric is a Yamabe soliton. Using these quantities we have obtained bound on the soliton constant. We have proved that the commutator of two soliton vector fields with the same metric in a given conformal class produces a Killing vector field. Also it is shown that the soliton vector field becomes a geodesic vector field if and only if the manifold is of constant curvature.
\end{abstract}
\maketitle
\section{Introduction}
\indent Let $M^n (n\ge 3)$ be a compact Riemannian manifold without boundary with a Riemannian metric $g$. The Yamabe flow, introduced by (\cite{HAMILTON1}, \cite{HAMILTON2}) is a way of deforming the initial metric $g$ in a given conformal class by means of the PDE
\begin{equation}\label{1.1}
\frac{\partial g}{\partial t}=-(R-r)g,
\end{equation}
where $R$ is the scalar curvature and $r$ is the average of $R$
\begin{equation*}
i.e \ \ \ r=\frac{\int_{M} R dv_g}{\int_{M} dv_g}.
\end{equation*}
Similar to the Ricci soliton, Yamabe solitons are special solution of (\ref{1.1}) moving purely by diffeomorphism and homothety and also arises as the blow-up limit of (\ref{1.1}). More precisely a complete Riemannian metric $g$ is called a Yamabe soliton if it satisfies
\begin{equation}\label{1.2}
\pounds_Vg_{ij}=2(-c+R)g_{ij},
\end{equation}
where $V$ is the soliton vector field and $c$ is the soliton constant.
If $c$ is a function in space, then $g$ is said to be almost Yamabe soliton.\\
The soliton is shrinking, steady and expanding according as $c<0,=0$ and $>0$ respectively. After the introduction in section $1$, section $2$ deals with some evolution equation of the Christoffel symbols, curvature tensor, Ricci tensor and scalar curvature. Section $3$ is devoted with the study commutator of two Yamabe solitons. Finally, we have obtained the condition of soliton vector field to be a geodesic vector field in the last section.
\section{Evolution of Geometric Quantities}
{\bf Lemma 2.1} If $(g,V,c)$ be a Yamabe soliton on a Riemannian manifold $M^n$, then
\begin{equation*}
(i) \  \pounds_V\Gamma^h_{ij}=(\nabla_jR)\delta^h_i+(\nabla_iR)\delta^h_j-\nabla^h R g_{ij},
\end{equation*}
\begin{equation*}
(ii) \ \pounds_VR^h_{kji}=(\nabla_k\nabla_iR)\delta^h_j-(\nabla_j\nabla_iR)\delta^h_k+(\nabla_j\nabla^hR)g_{ik}-(\nabla_k\nabla^hR)g_{ij},
\end{equation*}
\begin{equation*}
(iii) \ \pounds_VR_{ij}=g^{jk}\nabla_k\nabla_iR-\Delta R g_{ij}+\nabla_j\nabla^iR-n\nabla_j\nabla_iR,
\end{equation*}
\begin{equation*}
(iv) \ \pounds_VR=2(1-n)\Delta R,
\end{equation*}
where $R^h_{kji}$ and $R$ are respectively the curvature tensor, Ricci tensor and scalar curvature of $M$.\\
Proof: Using equation (\ref{1.2}) in the formula \cite{YANO}
\begin{equation*}
\pounds_V\Gamma^h_{ij}=\frac{1}{2}g^{hl}[\nabla_j(\pounds_Vg_{il})+\nabla_i(\pounds_Vg_{jl})-\nabla_l(\pounds_Vg_{ij})]
\end{equation*}
we have,
\begin{eqnarray}
\label{2.1}
 \pounds_V\Gamma^h_{ij} &=& g^{hl}[(\nabla_jR)g_{il}+(\nabla_iR)g_{jl}-(\nabla_lR)g_{ij}] \\
 \nonumber  &=&(\nabla_jR)\delta^h_i+(\nabla_iR)\delta^h_j-\nabla^h R g_{ij}.
\end{eqnarray}
This is the evolution of the Christoffel symbol or more preisely the connection.\\Having this hand and using commutation formula \cite{YANO}
\begin{equation*}
\nabla_k(\pounds_V\Gamma^h_{ij})-\nabla_j(\pounds_V\Gamma^h_{ik})=\pounds_VR^h_{kji}
\end{equation*}
we have
\begin{eqnarray}
\label{2.2}
\pounds_V R^h_{kji} &=&(\nabla_k\nabla_jR)\delta^h_i+(\nabla_k\nabla_iR)\delta^h_j-(\nabla_k\nabla^hR)g_{ij} \\
\nonumber &-&(\nabla_j\nabla_kR)\delta^h_i-(\nabla_j\nabla_iR)\delta^h_k+(\nabla_j\nabla^hR)g_{ik}\\
\nonumber &=& (\nabla_k\nabla_iR)\delta^h_j-(\nabla_j\nabla_iR)\delta^h_k+(\nabla_j\nabla^hR)g_{ik}-(\nabla_k\nabla^hR)g_{ij}.
\end{eqnarray}
Thus we get the evolution equation of curvature tensor along Yamabe soliton.\\
Again contracting (\ref{2.2}) we have,
\begin{equation}\label{2.3}
\pounds_VR_{ij}=g^{jk}\nabla_k\nabla_iR-\Delta R g_{ij}+\nabla_j\nabla^iR-n\nabla_j\nabla_iR.
\end{equation}
This gives the evolution equation for the Ricci tensor.\\
Again contracting (\ref{2.3}) we have,
\begin{eqnarray}\label{2.4}
\pounds_VR &=& \Delta R-n \Delta R+\Delta R-n\Delta R\\
\nonumber &=& 2(1-n)\Delta R.
\end{eqnarray}
This completes the proof of the lemma.\\
Now we can state our main theorem as the following:
\begin{theorem}
Let $M(g,V,c)$ be a compact Yamabe soliton. If the scalar curvature is bounded below by some positive constant $\alpha$, then we have the same lower bound for the soliton constant and the soliton is expanding.
\end{theorem}
\noindent {\bf Proof.} From (\ref{2.4}) we have,
\begin{equation}\label{2.5}
g(\nabla R,V)=(1-n)\Delta R.
\end{equation}
Since $M$ is compact, using divergence theorem\\
we have,
\begin{equation}\label{2.6}
\int_{M} g(\nabla R,V)dv_g = 0.
\end{equation}
Now note that,
\begin{eqnarray}
\label{6.a}
div(RV) &=& \nabla_i(RV^i)\\
\nonumber &=& g(\nabla R,V)+R \ div \ V.
\end{eqnarray}
and integrate (\ref{6.a}) over $M$ in order to get,
\begin{equation}\label{2.7}
\int_{M} [g(\nabla R,V)+R \ div \ V]dv_g=0.
\end{equation}
Now contracting equation (\ref{1.2}) we obtain
\begin{equation}\label{2.8}
div \ V=n(-c+R).
\end{equation}
Using (\ref{2.8}) and (\ref{2.6}) in (\ref{2.7}) we have,
\begin{equation*}
\int_{M} R \ div \ X \ dv_g=0,
\end{equation*}
\begin{equation*}
i.e \ \int_{M} R(-c+R) dv_g=0,
\end{equation*}
\begin{equation*}
i.e \ \ c=\frac{\int_{M} R^2 \ dv_g}{\int_{M} R \ dv_g}.
\end{equation*}
If $R\ge \alpha$ for some $\alpha>0$, then we have $c\ge \alpha$. This proves our theorem.
\section{Commutator of Soliton vector fields}
\indent We now consider two distinct Yamabe solitons with the same background Riemannian metric and consider commutator of there flow vector fields. More preceisely we prove that
\begin{theorem}
Let $M(g,V_1,c_1)$ and $M(g,V_2,c_2)$ be two distinct non-trivial Yamabe solitons. Then $[V_1,V_2]$ is a Killing vector field.
\end{theorem}
\noindent {\bf Proof.} Since $(g,V_1,c_1)$ and $(g,V_2,c_2)$ are Yamabe solitons, we have
\begin{equation}\label{3.1}
\begin{cases}
\pounds_{V_1}g= 2(R-c_1)g,& \mbox{ }\\
\pounds_{V_2}g= 2(R-c_2)g.& \mbox{}
\end{cases}
\end{equation}
Since the Yamabe solitons are unique, without loss of generality we may assume $c_1<c_2$. From (\ref{3.1}) we have,
\begin{equation*}
\pounds_{V_1-V_2}g=2(c_2-c_1)g.
\end{equation*}
This shows that $X = V_1-V_2$ is a homothetic vector field. Now,
\begin{eqnarray*}
\pounds_{[V_1,V_2]}g&=&\pounds_X\pounds_{V_2}g-\pounds_{V_2}\pounds_Xg\\
\nonumber &=& 2\pounds _X(R-c_2)g-2\pounds_{V_2}(c_2-c_1)g\\
\nonumber &=& 4(R-c_2)(c_2-c_1)g-4(c_2-c_1)(R-c_2)g\\
\nonumber &=& 0.
\end{eqnarray*}
i.e $[V_1,V_2]$ is Killing.\\
This proves the proposition.\\
\noindent {\bf Remark 1.} The theorem is also true for almost Yamabe solitons.
\section{Yamabe solitons and Geodesic vector fields}
\indent In this section we have studied solitons and obtained the condition for which its potential vector field becomes geodesic vector field in the following sense.
\begin{definition}
\cite{YN} A vector field $X$ on a Riemannian manifold $M$ is called a geodesic vector field if $\Box X=0$,
where $\Box$ is an operator acting on a smooth vector field $X^i$, given in local coordinates by
\begin{equation}\label{4.1}
\Box X^i=-(g^{jk}\nabla_j\nabla_kX^i+R^i_jX^j).
\end{equation}
\end{definition}
We now prove the following:
\begin{theorem}
If $(g,V,c)$ is a Yamabe soliton on $M$, then $V$ is a geodesic vector field if and only if $M$ is of constant scalar curvature.
\end{theorem}
\noindent {\bf Proof.} Tracing equation (\ref{1.2}) we have,
\begin{equation}\label{4.2}
\nabla_iV^i=(R-c)n.
\end{equation}
Taking covariant derivative of (\ref{4.2}) we have
\begin{equation}\label{4.3}
\nabla_j\nabla_iV^i=n\nabla_jR.
\end{equation}
Again writing (\ref{1.2}) as
\begin{equation*}
\nabla_jV^i+\nabla^iV_j=2(R-c)\delta^i_j
\end{equation*}
and taking covariant derivative it we get
\begin{equation}\label{4.4}
\nabla_i\nabla_jV^i+\nabla_i\nabla^iV_j=2\nabla_jR.
\end{equation}
Now subtracting (\ref{4.3}) from (\ref{4.4}) we have,
\begin{equation*}
\nabla_i\nabla_jV^i-\nabla_j\nabla^iV_i+\nabla_i\nabla^iV_j=(2-n)\nabla_jR.
\end{equation*}
Now applying Ricci formula we have
\begin{equation}\label{4.5}
\Box V=(n-2)\nabla_jR.
\end{equation}
This proves the theorem.\\
\noindent{\bf Remark 2.} As for a surface, the notion of Ricci soliton and Yamabe soliton coincides, it follows from (\ref{4.5}) that the soliton vector field is a geodesic vector field.

\vspace{0.5in}
\noindent Debabrata Chakraborty\\
Department of Mathematics, \\
Sidho Kanho Birsha University,\\
Purulia, 723104, West Bengal, India\\
E-mail: debabratamath@gmail.com

\vspace{0.1in}
\noindent Yadab Chandra Mandal\\
Department of Mathematics, \\
The University of Burdwan,\\
Burdwan, 713104, West Bengal, India \\
E-mail: myadab436@gmail.com

\vspace{0.1in}
\noindent Shyamal Kumar Hui\\
Department of Mathematics, \\
The University of Burdwan,\\
Burdwan, 713104, West Bengal, India \\
E-mail: skhui@math.buruniv.ac.in
\end{document}